\documentclass[a4paper,10pt]{article}
\usepackage{amsmath,amssymb}
\usepackage[utf8]{inputenc}
\usepackage[dvips]{epsfig}
\usepackage[dvips]{graphics}
\usepackage{color}

\newcommand{\dias}{\operatorname{Dias}}

\newcommand{\taq}{\mathcal{A}}
\newcommand{\TA}{\mathbb{A}}
\newcommand{\TAQ}{\overrightarrow{\mathbb{A}}}

\newcommand{\id}{\operatorname{Id}}

\newcommand{\bS}{\mathsf{S}}

\newtheorem{theorem}{Theorem}[section] 
\newtheorem{proposition}[theorem]{Proposition} 
 
\newtheorem{corollary}[theorem]{Corollary} 
\newtheorem{lemma}[theorem]{Lemma}

\newenvironment{proof}{\begin{trivlist}\item{\bf{Proof.}}}
  {\hfill\rule{2mm}{2mm}\end{trivlist}}

\title{The categorified Diassociative cooperad}
\author{F. Chapoton}
\date{\today}

\begin{document}

\maketitle

\begin{abstract}
  Using representations of quivers of type $\TA$, we define an
  anticyclic cooperad in the category of triangulated categories,
  which is a categorification of the linear dual of the Diassociative operad.
\end{abstract}

\section{Introduction}

The Diassociative operad has been introduced by Loday
\cite{loday_cras,loday_overview,loday_lnm}. It can be described as a
collection of free abelian groups $\dias(n)$ of rank $n$ and maps
$\circ_i$ from $\dias(m)\otimes \dias(n)$ to $\dias(m+n-1)$ satisfying
some kind of associativity. The composition maps $\circ_i$ have a
simple combinatorial description, using grafting of planar trees with
a distinguished path from the root to a leaf.

It has been shown in \cite{anticyclic} that one can endow this operad
with a refined structure of anticyclic operad. This means that there
exists a map of order $n+1$ on $\dias(n)$, with some compatibility
with the $\circ_i$ maps.

The aim of this article is to prove that this whole structure (or
rather its linear dual, which is an anticyclic cooperad) is the shadow of a
natural representation-theoretic object, related to the Dynkin
diagrams of type $\TA$.

We will not assume any knowledge of operads, but the interested reader
can consult \cite{ginzburg,loday_lnm,smirnov,stasheff} for basics of this
theory and \cite{markl,anticyclic} for the notion of anticyclic operad.

\smallskip

We first define a cooperad $\taq$ in the category of abelian categories.
This amounts to a collection of abelian categories $\taq_n$ for $n\geq 1$
and some functors $\nabla$ from these categories to products of two of
these categories. The categories $\taq_n$ involved are just the
categories of modules over the $\TAQ_n$ quivers. These are very
classical objects in representation theory.

The $\nabla$ functors are defined as tensor product by some specific
multiplicity-free bimodules. The axioms of cooperads are checked by
using a combinatorial description of the tensor product of such
bimodules.

At the level of the Grothendieck groups, one then checks that the
induced cooperad is the linear dual of the Diassociative operad. The
classes of simple modules correspond to the usual basis of $\dias(n)$
and the $\nabla$ functors give the linear dual of the $\circ_i$ maps.

\smallskip

As the $\nabla$ functors are given by the tensor product with
projective bimodules, they are exact. Going to the derived categories
$D\taq_n$, we prove that there is some compatibility between the
$\nabla$ functors and the Auslander-Reiten translations. At the level
of Grothendieck groups, this amounts to the structure of anticyclic
cooperad on the Diassociative cooperad.

\section{General facts}

\subsection{Quivers of type $\TA$}

For each integer $n\geq 1$, let $\TAQ_n$ be the quiver $1 \to 2 \to
\dots \to n$. This is a quiver on the graph of type $\TA_n$ in the
classification of Dynkin diagrams.

Let $k$ be a fixed ground field. Let $\taq_n$ be the category of
finite dimensional right-modules over the path algebra of $\TAQ_n$
over $k$. This is an abelian category, with a finite number of
isomorphism classes of indecomposable objects.

Let $D\taq_n$ be the bounded derived category of the category $\taq_n$. This
is a triangulated category, with a shift functor that will be denoted
by $S$. Indecomposable objects in $D\taq_n$ are just shifts of the
images of the indecomposable objects in $\taq_n$.

There exists a canonical self-equivalence of $D\taq_n$, called the
Auslander-Reiten translation and denoted by $\tau_n$.

The Nakayama functor $\nu_n$ is the composite $\tau_n S=S \tau_n$. This functor
maps, for each vertex $i$ of $\TAQ_n$, the projective module $P_i$ to
the injective module $I_i$.

\subsection{Products of quivers}

Let $\TAQ_{m_1,m_2,\dots,m_k}$ be the product quiver $\TAQ_{m_1}
\times \TAQ_{m_2} \times \dots \times \TAQ_{m_k}$. We consider this as
a quiver with relations by imposing all possible commutation
relations. 

A module over this quiver amounts to a module over the tensor product
of the path algebras of the quivers $\TAQ_{m_i}$. Therefore, one can
forget the action of some of the factors to define restricted modules.

As there is a canonical isomorphism of quivers from $\TAQ_{m,n}$ to
$\TAQ_{n,m}$, there are canonical equivalences between the
corresponding module and derived categories. Let us denote by $X$ these flip
functors. 

More generally, any permutation of the factors in a multiple product
of quivers $\TAQ_n$ give rise to a corresponding equivalence.

\subsection{Standard modules}

Let $M$ be a module over a quiver $\TAQ_{m_1,m_2,\dots,m_k}$. Let
$M_s$ be the vector space associated with a vertex $s$. One says that
$M$ is \textbf{multiplicity-free} if the dimension of $M_{s}$ is at
most $1$ for every vertex $s$. Let then $\bS(M)$ be the support of
$M$, which is the set of vertices $s$ such that $\dim M_s=1$.

Let $M$ be a multiplicity-free module. One says that $M$ is
\textbf{standard} if, for any two adjacent vertices $s,s'$ in
$\bS(M)$, the map between $M_s$ and $M_{s'}$ is an isomorphism.

To describe a standard module $M$ up to isomorphism, it is clearly
enough to give its support $\bS(M)$. The support of a standard module
cannot be arbitrary, because of the commuting conditions that must be
satisfied. One can then build back the module using copies of the
field $k$ and identity maps between them.


\subsection{Tensor product of projective standard modules}

There is a simple combinatorial description of the tensor product of
two projective standard modules.

Let us consider only the special case that we will
need. Let $M_{a\,;\,b,c}$ be a $\TAQ_a^{op} \times \TAQ_b \times
\TAQ_c$-module and let $M_{c\,;\,d,e}$ be a $\TAQ_c^{op} \times \TAQ_d
\times \TAQ_e$-module. Assume that $M_{a\,;\,b,c}$ is $\TAQ_{c}$-projective and
that $M_{c\,;\,d,e}$  is $\TAQ_c^{op}$-projective.

Then one can define the tensor product of $M_{a\,;\,b,c}$ and
$M_{c\,;\,d,e}$ over (the path algebra of) $\TAQ_c$. This is a
$\TAQ_a^{op} \times \TAQ_b \times \TAQ_d \times \TAQ_e$-module denoted
by $M_{a\,;\,b,c} \otimes_{\TAQ_c} M_{c\,;\,d,e}$.

Assume that $M_{a\,;\,b,c}$ and $M_{c\,;\,d,e}$ are standard modules
with support $\bS_{a\,;\,b,c}$ and $\bS_{c\,;\,d,e}$. Let us define a
set $\bS_{a\,;\,b,c}\times_c\bS_{c\,;\,d,e}$ as follows:
\begin{equation*}
  \bS_{a\,;\,b,c}\times_c\bS_{c\,;\,d,e}=\{
(\alpha,\beta,\delta,\epsilon)\mid \exists \gamma \, (\alpha,\beta,\gamma)\in\bS_{a\,;\,b,c} \text{ and }
(\gamma,\delta,\epsilon)\in\bS_{c\,;\,d,e}
\}.
\end{equation*}

\begin{proposition}
  \label{def_tensor}
  The tensor product $M_{a\,;\,b,c} \otimes_{\TAQ_c} M_{c\,;\,d,e}$ is
  isomorphic to the standard module with support
  $\bS_{a\,;\,b,c}\times_c\bS_{c\,;\,d,e}$.
\end{proposition}
\begin{proof}
  The tensor product over the field $k$ has a basis indexed by tuples
  $(\alpha,\beta,\gamma,\gamma',\delta,\epsilon)$ with
  $(\alpha,\beta,\gamma)\in\bS_{a\,;\,b,c}$ and
  $(\gamma',\delta,\epsilon)\in\bS_{c\,;\,d,e}$. Then one has to take
  the quotient by the action of the idempotents and arrows of the
  quiver $\TAQ_c$. Abusing notation, we will identify tuples with the
  corresponding vectors.

  By the action of the idempotents in the path algebra, one can see
  that all vectors $(\alpha,\beta,\gamma,\gamma',\delta,\epsilon)$
  with $\gamma\not=\gamma'$ vanish in the tensor product.

  There remains to quotient by the action of the arrows. This means
  that one has to identify
  $(\alpha,\beta,\gamma,\gamma,\delta,\epsilon)$ and
  $(\alpha,\beta,\gamma+1,\gamma+1,\delta,\epsilon)$, provided that
  one has $(\alpha,\beta,\gamma)\in \bS_{a\,;\,b,c}$ and
  $(\gamma+1,\delta,\epsilon)\in \bS_{c\,;\,d,e}$.

  By the hypothesis of projectivity, in this situation, one also has
  $(\alpha,\beta,\gamma+1)\in \bS_{a\,;\,b,c}$ and
  $(\gamma,\delta,\epsilon)\in \bS_{c\,;\,d,e}$. Hence both
  $(\alpha,\beta,\gamma,\gamma,\delta,\epsilon)$ and
  $(\alpha,\beta,\gamma+1,\gamma+1,\delta,\epsilon)$ are non-zero vectors.

  Therefore, the tensor product has a basis indexed by tuples
  $(\alpha,\beta,\delta,\epsilon)$ such that there exists $\gamma$
  with $(\alpha,\beta,\gamma)\in\bS_{a\,;\,b,c}$ and
  $(\gamma,\delta,\epsilon)\in\bS_{c\,;\,d,e}$.

  One can also see by construction that indeed the module is standard.
\end{proof}

\subsection{Fiber-reversal and action of $\tau$}

Let $N_n$ be the standard $\TAQ_n^{op} \times \TAQ_n$ module
with support 
\begin{equation}
  \{(i,j)\in[1,n]\times[1,n] \mid i \geq j\}.
\end{equation}
Note that $N_n$ is injective as a $\TAQ_n^{op}$-module and as a $\TAQ_n$-module.

\begin{lemma}
  The Nakayama functor $\nu_n$ on the category $D\taq_n$ is the derived
  tensor product $? \otimes^{L}_{\TAQ_n} N_n$.
\end{lemma}
\begin{proof}
  This follows from the fact that the image by $\nu$ of the projective
  module $P_i$ is the injective module $I_i$, by the standard way of
  representing functors by bimodules.
\end{proof}

\begin{figure}
  \begin{center}
    \epsfig{file=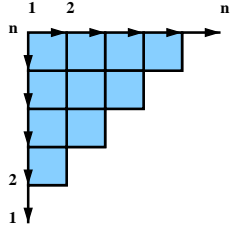,height=3cm} 
    \caption{The bimodule $N_{6}$ corresponding to Nakayama functor $\nu_6$}
    \label{exemple_nakayama}
  \end{center}
\end{figure}

\medskip

Let us now introduce some operations on support sets.

Let $\bS$ be a subset in the product
$[1,m_1]\times\dots\times[1,m_k]$. Fix an index $i$. Assume that $\bS$
is projective in the direction $i$, \textit{i.e.} that
 \begin{equation}
  \text{if } (j_1,\dots,j_i,\dots,j_k)\in\bS \text{ then }
 (j_1,\dots,\ell,\dots,j_k)\in\bS \text{ for all }\ell\geq j_i.
 \end{equation}

 The \textbf{fiber-reversal} of $\bS$ in the direction $i$ is
\begin{equation}
 \{(j_i,\dots,j_i,\dots,j_k)\in
 [1,m_1]\times\dots\times[1,m_k] \mid
 (j_1,\dots,j_i-1,\dots,j_k)\not\in\bS \}.
\end{equation}

Note that the fiber-reversal of $\bS$ in direction $i$ is never
disjoint from $\bS$, and really depends on the index $i$.

One can give a similar definition of the fiber-reversal if the set
$\bS$ is injective in the direction $i$, \textit{i.e.} if the
following condition holds:
\begin{equation}
  \text{if } (j_1,\dots,j_i,\dots,j_k)\in\bS \text{ then }
  (j_1,\dots,\ell,\dots,j_k)\in\bS \text{ for all }\ell\leq j_i.
\end{equation}

Let us now describe the (derived) tensor product with $N_n$. We
consider only the special case that we will need.

Let $M_{n;c;d}$ be a $\TAQ_n^{op}\times \TAQ_{c} \times \TAQ_d$
standard module with support $\bS_{n;c;d}$. Assume that $M_{n;c,d}$ is
projective as a $\TAQ_n^{op}$-module.
\begin{proposition}
  The derived tensor product of $N_n \otimes_{\TAQ_n}^L M_{n;c,d}$ is
  isomorphic to the standard module with support the fiber-reversal
  of $\bS_{n;c;d}$ in the direction of length $n$.
\end{proposition}

\begin{proof}
  The tensor product $N_n \otimes_{k} M_{n;c,d}$ has a basis indexed
  by tuples $(\alpha,\beta,\beta',\gamma,\delta)$ with
  $\alpha \geq \beta$ and $(\beta',\gamma,\delta)\in\bS_{n;c;d}$. Abusing
  notation, we will identify tuples with the corresponding vectors.

  Using the idempotents in the path algebra, the tensor product over
  $\TAQ_n$ is spanned by tuples $(\alpha,\beta,\beta,\gamma,\delta)$.
  Then one has to identify $(\alpha,\beta,\beta,\gamma,\delta)$ and
  $(\alpha,\beta+1,\beta+1,\gamma,\delta)$ as soon as $\alpha \geq
  \beta$ and $(\beta+1,\gamma,\delta)\in \bS_{n;c,d}$.

  Using the hypothesis that $M_{n;c,d}$ is projective, one has in this
  situation that $(\beta,\gamma,\delta)\in \bS_{n;c,d}$. 

  The only case where one of the two vectors
  $(\alpha,\beta,\beta,\gamma,\delta)$ and
  $(\alpha,\beta+1,\beta+1,\gamma,\delta)$ is zero and the other is
  not zero happens if $\beta+1>\alpha$, \textit{i.e.} $\alpha=\beta$.

  It follows that the vector $(\alpha,\beta,\beta,\gamma,\delta)$ is
  mapped to zero in the tensor product over $\TAQ_n$ if and only if
  $(\alpha+1,\gamma,\delta) \in \bS(M)$ and are otherwise just
  identified. This is exactly the definition of the fiber-reversal of
  $\bS_{n;c,d}$ in the first direction.

  One can easily check that the tensor product is standard.
\end{proof}

\section{The $\nabla$ functors on module categories}

In this section, we define a cooperad structure on the collection of
module categories $(\taq_n)_{n\geq 1}$. This means that several
functors $\nabla$ are introduced and some kind of associativity
properties are proved.

Let $n\geq 1$ be an integer and let $m,i$ be integers such that $1
\leq i \leq m$.

Consider the set $\bS_{m;i}^n$ of integer triples
$(\gamma,\mu,\nu)$ in $[1,m+n-1]\times[1,m]\times[1,n]$ such that
\begin{align*}
 &( \mu \leq i-1 \text{ and } \gamma \leq \mu) \\
  \text{ or } &(\mu =i \text{ and } \gamma \leq i+\nu-1)\\
  \text{ or } &(i+1 \leq \mu  \text{ and } \gamma \leq \mu+n-1).
\end{align*}

\begin{figure}
  \begin{center}
    \scalebox{1.2}{\input{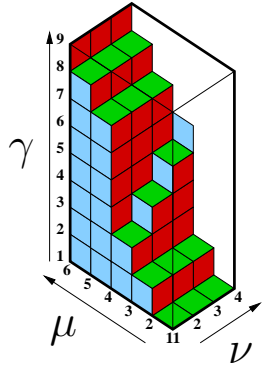}}
    \caption{The module $M_{6;3}^{4}$ and its symbolic description}
    \label{exemple_module}
  \end{center}
\end{figure}

This is illustrated in Figure \ref{exemple_module} with $m=6,n=4$ and $i=3$.

For later use, here is an equivalent description of $\bS_{m;i}^{n}$ :
\begin{align}
  \label{description_alt}
 &( \gamma \leq i \text{ and } \gamma \leq \mu) \\
  \text{ or } &(i+1 \leq \gamma\leq i+n-1 \text{ and } \gamma \leq
  i+\nu-1)\text{ and } i\leq \mu\\
  \text{ or } &(i+1 \leq \gamma\leq i+n-1 \text{ and }  i+\nu \leq
  \gamma )\text{ and } i+1\leq \mu\\
  \text{ or } &(i+n \leq \gamma  \text{ and } \gamma-n+1 \leq \mu).
\end{align}

One can easily check that the set $\bS_{m;i}^n$ has the following
property : if $(\gamma,\mu,\nu)\in \bS_{m;i}^n$ and if
$(\gamma',\mu',\nu')\in [1,m+n-1]\times[1,m]\times[1,n]$ is such that
$\gamma' \leq \gamma$, $\mu \leq \mu'$ and $\nu \leq \nu'$, then
$(\gamma',\mu',\nu',\gamma')\in \bS_{m;i}^n$.

This implies that one can define a representation $M_{m;i}^n$ of the
quiver $\TAQ_{m+n-1}^{op} \times \TAQ_m \times \TAQ_n$ as the standard
module with support $\bS_{m;i}^n$. Note that $M_{m;i}^n$ is projective
with respect to $\TAQ_m$, $\TAQ_n$ and $\TAQ_{m+n-1}^{op}$.

Let then $\nabla_{m;i}^n$ be the functor from $\taq_{m+n-1}$ to
$\taq_{m,n}$ defined by the tensor product over
$\TAQ_{m+n-1}$ by $M_{m;i}^{n}$ :
\begin{equation}
  \nabla_{m;i}^n = ? \otimes_{\TAQ_{m+n-1}} M_{m;i}^{n}.
\end{equation}

Note that $\nabla_{1;1}^n$ is the identity functor of $\taq_n$ and that $\nabla_{m;i}^{1}$ is the identity functor of $\taq_m$ for all $i$.

\subsection{Relation to the $\dias$ cooperad}

Let $n$ be an integer and let $1\leq j \leq n$. Let $S_j^n$ be the
simple $\TAQ_{n}$-module supported at vertex $j$. Let $P_j^n$ be the
projective $\TAQ_{n}$-module associated with vertex $j$. The class of
a module $M$ in the Grothendieck group $K^0(\taq_n)$ of $\taq_n$ will
be denoted by $[M]$. The elements $[S_j^n]$ for $1\leq j \leq n$ form
a basis of $K^0(\taq_n)$.

Let us now compute the class $[\nabla_{m;i}^{n}(S_j^{m+n-1})]$.

From the explicit description of the module $M_{m;i}^{n}$, one has
\begin{equation}
  [\nabla_{m;i}^n (P_j^{m+n-1})]=
  \begin{cases}
    [P_{j,1}^{m,n}] &\text{ if } 1\leq j \leq i,\\
    [P_{i+1,1}^{m,n}]+[P_{i,j-i+1}^{m,n}]-[P_{i+1,j-i+1}^{m,n}] 
                              &\text{ if } i+1\leq j\leq i+n-1,\\
    [P_{j-n+1,1}^{m,n}] &\text{ if } i+n \leq j \leq m+n-1,
  \end{cases}
\end{equation}
where $P_{i,j}^{m,n}$ is the projective module associated with vertex
$(i,j)$ of $\TAQ_m \times \TAQ_n$.

Using a projective resolution of the simple modules, one deduces that
\begin{equation}
   [\nabla_{m;i}^n (S_j^{m+n-1})]=
  \begin{cases}
   \sum_{k=1}^{n} [S_{j,k}^{m,n}]&\text{ if } 1\leq j \leq i-1,\\
   [S_{i,j-i+1}^{m,n}]&\text{ if } i\leq j\leq i+n-1,\\
   \sum_{k=1}^{n} [S_{j-n+1,k}^{m,n}]&\text{ if } i+n \leq j \leq m+n-1,
  \end{cases} 
\end{equation}
where $S^{m,n}_{i,j}$ is the simple module at vertex $(i,j)$ for
$\TAQ_{m} \times \TAQ_{n}$.

Taking the linear dual basis $e$ of the basis $[S]$, one finds that
the linear dual maps $\circ$ to the $\nabla$ maps are given by
\begin{equation}
  \circ_{m;i}^{n}(e^m_j \otimes e^n_k) =
  \begin{cases}
    e^{m+n-1}_{j}&\text{ if }i>j,\\
    e^{m+n-1}_{i+k-1}&\text{ if }i=j,\\
    e^{m+n-1}_{j+n-1}&\text{ if }i<j.
  \end{cases}
\end{equation}

This is exactly the usual description of the Diassociative operad, in
the usual basis $e$ of $\dias(n)$, see \cite[\S 3]{anticyclic}.

\section{Cooperadic properties of $\nabla$ functors}

One has to check two different axioms to prove that the $\nabla$
functors define a cooperad. Let us call them the commutativity axiom
and the associativity axiom.

\subsection{Commutativity axiom}

Let $m,n,p$ and $i,j$ be integers such that $1\leq i <j \leq m$.

\begin{proposition}
  The modules $M$ have the following property : there is an
  isomorphism
   \begin{equation}
    M_{m+p-1;i}^n   \otimes_{\TAQ_{m+p-1}} M_{m;j}^p\simeq 
     M_{m+n-1;j+n-1}^p \otimes_{\TAQ_{m+n-1}}M_{m;i}^n ,
   \end{equation}
   where both sides are $\TAQ_{m+n+p-2}^{op}\times\TAQ_m\times \TAQ_n
   \times \TAQ_p$-modules.
\end{proposition}

\begin{proof}
  As the modules $M$ are standard and projective, their tensor
  products can be described using their supports. According to the
  description of tensor products in Prop. \ref{def_tensor}, one
  therefore has to compute and compare the sets $\bS_{m;j}^p
  \times_{{m+p-1}} S_{m+p-1;i}^n$ and $\bS_{m;i}^n \times_{{m+n-1}}
  S_{m+n-1;j+n-1}^p$.

  By an elementary computation with boolean combinations of
  inequalities, one can show that both sides are given by the set of
  $(\gamma,\mu,\nu,\pi)$ in
  $[1,m+n+p-2]\times[1,m]\times[1,n]\times[1,p]$ such that
  \begin{align*}
    &( \mu \leq i-1 \text{ and } \gamma \leq \mu) \\
    \text{ or } &(\mu =i \text{ and } \gamma \leq i+\nu-1)\\
    \text{ or } &(i+1 \leq \mu \leq j-1 \text{ and } \gamma \leq
    \mu+n-1)\\
    \text{ or } &(\mu =j \text{ and } \gamma \leq j+n+\pi-2)\\
    \text{ or } &(j+1 \leq \mu \text{ and } \gamma \leq
    \mu+n+p-2).
  \end{align*}
\end{proof}

\begin{corollary}
 The functors $\nabla$ have the following property : there is a
 natural transformation
 \begin{equation}
  (\id_m \times X)(\nabla_{m;j}^p \times \id_n)\nabla_{m+p-1;i}^n \simeq (\nabla_{m;i}^n \times \id_p) \nabla_{m+n-1;j+n-1}^p.
 \end{equation}
\end{corollary}

\subsection{Associativity axiom}

Let $m,n,p$ and $i,j$ be integers such that $1\leq i \leq m$ and $1 \leq j \leq n$.

\begin{proposition}
  The modules $M$ have the following property : there is an
  isomorphism
  \begin{equation}
   M_{m;i}^{n+p-1}  \otimes_{\TAQ_{n+p-1}} M_{n;j}^{p} \simeq  M_{m+n-1;j+i-1}^p\otimes_{\TAQ_{m+n-1}} M_{m;i}^n ,
  \end{equation}
  where both sides are $\TAQ_{m+n+p-2}^{op}\times\TAQ_m\times \TAQ_n
  \times \TAQ_p$-modules.
\end{proposition}

\begin{proof}
  As in the previous section, one just has to compute the supports of these modules. One can check that they both give the
  set of $(\mu,\nu,\pi,\gamma)$ in $[1,m]\times[1,n]\times[1,p]\times[1,m+n+p-2]$ such that
  \begin{align*}
    &( \mu \leq i-1 \text{ and } \gamma \leq \mu) \\
    \text{ or } &(\mu =i \text{ and } \nu\leq j-1 \text{ and }\gamma \leq i+\nu-1)\\
    \text{ or } &(\mu =i \text{ and }\nu=j \text{ and } \gamma \leq
    i+j+\pi-2)\\
    \text{ or } &(\mu =i \text{ and }j+1 \leq \nu \text{ and } \gamma \leq i+\nu+p-2)\\
    \text{ or } &(i+1 \leq \mu \text{ and } \gamma \leq
    \mu+n+p-2).
  \end{align*}
\end{proof}

\begin{corollary}
 The functors $\nabla$ have the following property : there is a
 natural transformation
 \begin{equation}
  (\id_m \times \nabla_{n;j}^{p})\nabla_{m;i}^{n+p-1} \simeq (\nabla_{m;i}^n \times \id_p) \nabla_{m+n-1;j+i-1}^p.
 \end{equation}
\end{corollary}

\section{Relations between $\nabla$ and $\tau$}

Let us consider the functors $\nabla_{m;i}^{n}$ as the derived tensor
product with $M_{m;i}^{n}$. As the modules $M_{m;i}^{n}$ are
projective in every direction, this is just the usual tensor product.
Therefore, we obtain a cooperad structure on the collection of derived
categories $(D\taq_n)_{n\geq 1}$.

In this section, we define an anticyclic cooperad structure on the
collection of derived categories $(D\taq_n)_{n\geq 1}$. This means
that some compatibility properties hold between the functors $\nabla$
and the Auslander-Reiten translations $\tau$. We will rather work with
the Nakayama functors $\nu$, described as derived tensor product with
the modules $N$.

There are two different axioms for the notion of anticyclic operad :
let us call them the border axiom and the inner axiom.

\subsection{Border axiom}

\begin{figure}
  \begin{center}
    \scalebox{1}{\input{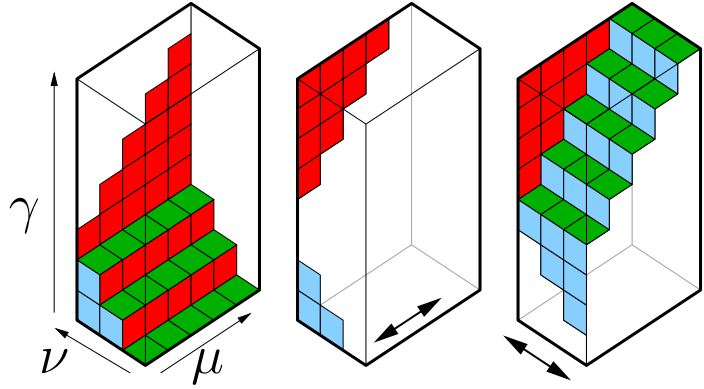}}
    \caption{The module $M_{4;4}^{6}$, its fiber-reversal in the
      direction of length $6$ and the fiber-reversal of the result
      in the direction of length $4$}
    \label{exemple_bord1}
  \end{center}
\end{figure}

\begin{proposition}
  The fiber-reversal of $\bS_{m;1}^{n}$ in the direction of length
  $m+n-1$ is equal to the fiber-reversal in the direction of length
  $n$ of the fiber-reversal in the direction of length $m$ of
  $\bS_{n;n}^{m}$. In terms of modules, this means that
  \begin{equation}
    N_{m+n-1} \otimes^L_{\TAQ_{m+n-1}} M_{m;1}^n \simeq (M_{n;n}^m
    \otimes^L_{\TAQ_m} N_m) \otimes^L_{\TAQ_n} N_n.
  \end{equation}
\end{proposition}

\begin{proof}
  Let us first compute the fiber-reversal of  $\bS_{m;1}^{n}$ in the
  direction of $\gamma$ of length
  $m+n-1$ . One easily gets
  \begin{align*}
   &(\mu =1 \text{ and } \gamma \geq \nu)\\
  \text{ or } & \gamma \geq \mu+n-1.
\end{align*}

Let us then compute the fiber-reversal of $\bS_{n;n}^{m}$ in the
direction of $\mu$ of length $m$. One gets
\begin{align*}
 &( \mu = 1 \text{ and } \gamma \leq \nu) \\
  \text{ or } &(\nu =n \text{ and } \gamma \geq n+\mu-1).
\end{align*}

Then one can compute the fiber-reversal of this set in the
direction of $\nu$ and check the expected result.
\end{proof}

\begin{figure}
  \begin{center}
    \scalebox{1}{\input{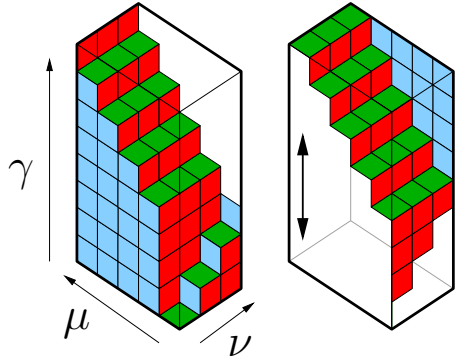}}
    \caption{The module $M_{6;1}^{4}$, its top-bottom fiber-reversal}
    \label{exemple_bord2}
  \end{center}
\end{figure}

\begin{corollary}
  The functors $\nabla$ satisfy
 \begin{equation}
   \label{anticyclic_reversal}
  \nabla_{m;1}^{n} \tau_{m+n-1} \simeq S X (\tau_n \times \tau_m) \nabla_{n;n}^{m}.
 \end{equation}
\end{corollary}

\subsection{Inner axiom}

Let us assume here that $2 \leq i \leq m$.

\begin{figure}
  \begin{center}
    \scalebox{1}{\input{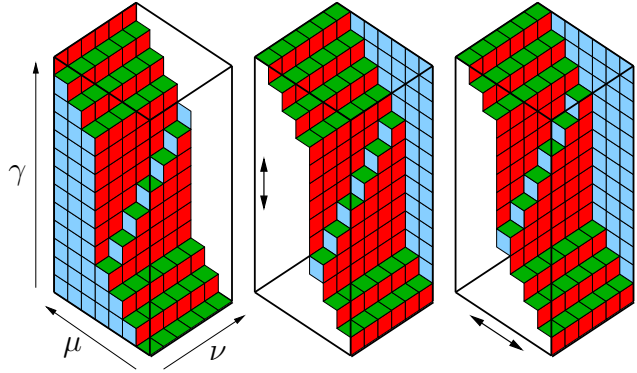}}
    \caption{The module $M_{8;4}^{7}$, its top-bottom fiber-reversal and its
      fiber-reversal in the direction of length $8$}
    \label{exemple_interieur}
  \end{center}
\end{figure}

\begin{proposition}
  The fiber-reversal of $\bS_{m;i}^{n}$ in the direction of length
  $m+n-1$ coincides with the fiber-reversal of $\bS_{m;i-1}^{n}$ in
  the direction of length $m$. In terms of modules, this means
  \begin{equation}
    N_{m+n-1} \otimes^L_{\TAQ_{m+n-1}} M_{m;i}^n \simeq M_{m;i-1}^n
    \otimes^L_{\TAQ_m} N_m.
  \end{equation}
\end{proposition}

\begin{proof}
  On the one hand, for the fiber-reversal of $\bS_{m;i}^{n}$ in the
  direction $\gamma$ of length $m+n-1$, one easily gets
  \begin{align*}
     &( \mu \leq i-1 \text{ and } \gamma \geq \mu) \\
  \text{ or } &(\mu =i \text{ and } \gamma \geq i+\nu-1)\\
  \text{ or } &(i+1 \leq \mu  \text{ and } \gamma \geq \mu+n-1).
  \end{align*}

  On the other hand, using the alternative description
  (\ref{description_alt}) of $S_{m;i-1}^n$, the fiber-reversal of
  $\bS_{m;i-1}^{n}$ in the direction $\mu$ of length $m$ is
  \begin{align*}
  \label{description_alt}
 &( \gamma \leq i-1 \text{ and } \gamma \geq \mu) \\
  \text{ or } &(i \leq \gamma\leq i+n-2 \text{ and } \gamma \leq
  i+\nu-2)\text{ and } i-1\geq \mu\\
  \text{ or } &(i \leq \gamma\leq i+n-2 \text{ and }  i-1+\nu \leq
  \gamma )\text{ and } i\geq \mu\\
  \text{ or } &(i+n-1 \leq \gamma  \text{ and } \gamma-n+1 \geq \mu).
  \end{align*}
  It is then a routine matter to check that they are indeed equal.  
\end{proof}

\begin{corollary}
  The functors $\nabla$ satisfy
 \begin{equation}
  \label{anticyclic_shift}
 \nabla_{m;i}^{n} \tau_{m+n-1} \simeq (\tau_m \times \id) \nabla_{m;i-1}^{n},
\end{equation}
for $2 \leq i \leq m$.
\end{corollary}

\bibliographystyle{alpha}
\bibliography{catoperad}

Remerciements : Merci {\`a} Bernhard Keller pour ses conseils. Ce travail
a {\'e}t{\'e} partiellement financ{\'e} par un programme blanc de l'ANR.

\begin{center}
Fr{\'e}d{\'e}ric Chapoton\\
Universit{\'e} de Lyon ; Universit{\'e} Lyon 1 ;\\
CNRS, UMR5208, Institut Camille Jordan,\\
43 boulevard du 11 novembre 1918,\\
F-69622 Villeurbanne-Cedex, France
\end{center}

\end{document}